\documentclass[11pt,reqno]{amsart}
\usepackage{indentfirst, amssymb,amsmath,amsthm}
\newtheorem{theo}{Theorem}[section]
\newtheorem{exm}{Example}[section]
\newtheorem{lem}{Lemma}

\newtheorem{res}{Result}[section]
\newtheorem{defi}{Definition}
\newtheorem{obs}{Observation}

\newtheorem{nota}{Notation}

\newcommand{\be}{\begin{equation}}
\newcommand{\ee}{\end{equation}}
\newcommand{\beas}{\begin{eqnarray*}}
\newcommand{\eeas}{\end{eqnarray*}}
\newcommand{\bea}{\begin{eqnarray}}
\newcommand{\eea}{\end{eqnarray}}

\numberwithin{equation}{section}

\begin{document}

\setlength{\unitlength}{1mm} \baselineskip .52cm
\setcounter{page}{1}
\pagenumbering{arabic}
\title[On Lebesgue Property for Fuzzy Metric Spaces]{On Lebesgue Property for Fuzzy Metric Spaces}

\author[Sugata Adhya and A. Deb Ray]{Sugata Adhya and A. Deb Ray}

\address{Department of Mathematics, The Bhawanipur Education Society College. 5, Lala Lajpat Rai Sarani, Kolkata 700020, West Bengal, India.}
\email {sugataadhya@yahoo.com}

\address{Department of Pure Mathematics, University of Calcutta. 35, Ballygunge Circular Road, Kolkata 700019, West Bengal, India.}
\email {debrayatasi@gmail.com}

\maketitle

\begin{abstract}
    We provide several characterizations of the Lebesgue property for fuzzy metric spaces. It is known that a fuzzy metric space is Lebesgue if and only if every real-valued continuous function is uniformly continuous. Here we show that it suffices to examine uniform continuity of bounded real-valued continuous functions for characterizing Lebesgue property in fuzzy setting.\\
\end{abstract}


\section{\textbf{Introduction}}

For a long period of time, several definitions of fuzzy metric spaces were proposed towards establishing a consistent theory of metric fuzziness (\cite{1}, \cite{2}, \cite{3}). In particular, George and Veeramani made an appealing modification (\cite{ver1}, \cite{ver2}) to a definition proposed in \cite{km}. Using this modified version, they proved that every fuzzy metric induces a Hausdroff, first countable topology, which was later realized to be metrizable \cite{uni}. On the other hand, it was shown that if this fuzzy metric is complete, then the induced topology is completely metrizable \cite{uni}. In this context, the modified definition turned out to be an appropriate notion for metric fuzziness.

Several topological and metric concepts were extended in this new structure (\textit{see} e.g. \cite{ver1}, \cite{ver2}, \cite{uni}). In particular, the notion of uniform continuity was introduced by George and Veeramani in \cite{cont} where it was shown that every real-valued continuous function on a compact fuzzy metric space is uniformly continuous. Since the converse is not necessarily true, it became a point of interest to characterize those fuzzy metric spaces on which real-valued continuous functions are uniformly continuous. In \cite{uc}, Gregori, Romaguera and Sapena proposed the notion of Lebesgue property for metric fuzziness. They not only gave satisfactory extensions of many important properties of Lebesgue metric spaces in this new structure, but also obtained the solution of the above problem by concluding that the class of fuzzy metric spaces having Lebesgue property agrees with the class of fuzzy metric spaces on which real-valued continuous functions are uniformly continuous. 

In this paper, we extend the above study by obtaining several important characterizations of Lebesgue property for fuzzy metric spaces. In particular, we have shown that it is sufficient to examine the uniform continuity of bounded real-valued continuous functions to characterize Lebesgue property and have established that the class of fuzzy metric spaces on which real-valued continuous functions are uniformly continuous coincides with the class of fuzzy metric spaces on which bounded real-valued continuous functions are uniformly continuous.

\section{\textbf{Preliminaries}}

This section is aimed at collecting necessary definitions and facts on fuzzy metric spaces, to be required subsequently. For notions concerning general topology, we refer the reader to \cite{will}.

\begin{defi}
\normalfont \cite{tn} A binary operation $*$ on $[0,1]$ is called a continuous $t$-norm, if 

a) $*$ is associative and commutative,

b) $*$ is continuous 

c) $\forall~a\in[0,1],~a*1=a,$

d) $\forall~a,b,c,d\in[0,1],~a\le b,~c\le d\implies a*c\le b*d.$
\end{defi}

\begin{defi}
\normalfont \cite{ver1,ver2} A fuzzy metric space is an ordered
triple $(X,M,*)$ where $X$ is a nonempty set, $*$ is a continuous $t$-norm and $M$ is a fuzzy set of
$X\times X\times(0,\infty)$ such that, for all $x,y,z\in X$ and $s,t>0,$ the following conditions hold:

a) $M(x,y,t)>0,$

b) $M(x,y,t)=1\iff x=y,$

c) $M(x,y,t)=M(y,x,t)$

d) $M(x,y,t)*M(y,z,s)\le M(x,z,t+s)$

e) $M(x,y,.):(0,\infty)\to[0,1]$ is continuous.

In this case, $(M,*)$ is said to be a fuzzy metric on $X.$
\end{defi}

\begin{lem}
\normalfont \cite{tn} For all $x,y\in X,$ $M(x,y,.):(0,\infty)\to[0,1]$ is nondecresing. 
\end{lem}

Let $(X,d)$ be a metric space. Define $M_d:X\times X\times(0,\infty)\to[0,1]$ by $$M_d(x,y,t)=\frac{t}{t+d(x,y)}.$$ Then $(X,M_d,*)$ forms a fuzzy metric space, $*$ being the usual multiplication of real numbers. (\cite{ver1})

\begin{defi}
\normalfont $(M_d,*)$ is said to be the fuzzy metric induced by $d$ \cite{ver1} and $(X,M_d,*),$ to be the induced fuzzy metric space.
\end{defi}

Let $(X,M,*)$ be a fuzzy metric space. For $x\in X,r\in(0,1)$ and $t>0$ denote $B(x,r,t)=\{y\in X:M(x,y,t)>1-r\}.$ Then $\{B(x,r,t):x\in X,r\in(0,1),t>0\}$ forms a base for some topology $\tau_M$ on $X$ \cite{ver1}.

\begin{defi}
\normalfont $\tau_M$ is called the topology induced by $(M,*).$
\end{defi}

\begin{res}\label{res1}
\normalfont \cite{ver1} Let $(X,d)$ be a metric space. Then $\tau_{M_d}$ coincide with the topology induced by the metric $d$ on $X.$
\end{res}

In view of Result \ref{res1} it is clear, if $(X,\tau)$ is a metrizable topological space then there is a fuzzy metric $(M,*)$ on $X$ such that $\tau_M=\tau.$

\begin{defi}
\normalfont \cite{sub} Let $(X, M, *)$ be a fuzzy metric space and $A\subset X.$ If $M_A$ denotes the restriction of $M$ on $A\times A\times(0,\infty)$ then $(M_A,*)$ defines a fuzzy metric on $A.$ The associated fuzzy metric space $(A, M_A, *)$ is called the fuzzy metric subspace of $(X, M, *)$ on $A.$ 
\end{defi}

It is immediate that $\tau_{M_A}=(\tau_M)_A,$ $(\tau_M)_A$ being the subspace topology on $A$ induced by $\tau_M.$

It has been shown in \cite{ver1} that given a fuzzy metric space $(X,M,*),$ $(X,\tau_M)$ is a Hausdroff, first countable topological space. 

In \cite{uni}, it has been shown that $\{U_n:n\in\mathbb N\}$ forms a base for some uniformity $\mathcal{U}_M$ on $X$ such that $\mathcal{U}_M$ is compatible with $\tau_M,$ where $U_n=\left\{(x,y):M\left(x,y,\frac{1}{n}\right)>1-\frac{1}{n}\right\},~\forall~n\in\mathbb N,$ and thereby, much stronger result has been obtained as follows:

\begin{theo}
\normalfont \cite{uni} Let $(X,M,*)$ be a fuzzy metric space. Then $(X,\tau_M)$ is metrizable.
\end{theo}

\begin{defi}
\normalfont Let $(X,M,*),~(Y,N,\star)$ be fuzzy metric spaces. 

(a) A mapping $f:X\to Y$ is called uniformly continuous \cite{cont} if for $\epsilon\in(0,1),\delta>0$ there exist $r\in(0,1),s>0$ such that $\forall~x,y\in X,~M(x,y,s)>1-r\implies~N(f(x),f(y),\delta)>1-\epsilon.$

(b) A mapping $f:X\to\mathbb R$ is called $\mathbb R$-uniformly continuous \cite{uc} if for $\epsilon>0$ there exist $r\in(0,1),s>0$ such that $\forall~x,y\in X,~M(x,y,s)>1-r\implies~|f(x)-f(y)|<\epsilon.$
\end{defi}

Clearly $f:(X,M,*)\to (Y,N,\star)$ is uniformly continuous $\implies f:(X,\tau_M)\to(Y,\tau_N)$ is continuous.

\begin{defi}
\normalfont \cite{fd} Let $(X,M,*)$ be a fuzzy metric space. A subset $A$ of $X$ is called fuzzy uniformly discrete in $X$ if there exist $r\in(0,1),s>0$ such that $M(x,y,s)<1-r,~\forall~x,y~(x\ne y)\in A.$
\end{defi}

One may observe the following:

\begin{obs}
\normalfont Given a fuzzy metric space $(X,M,*)$, a subset $A$ of $X$ is fuzzy uniformly discrete if and only if there exists $r\in(0,1),s>0$ such that $B(a,r,s)\cap A=\{a\}~\forall~a\in A.$
\end{obs}

A sequence $(x_n)$ in a fuzzy metric space $(X,M,*)$ converges to $x\in X$ [\textit{briefly,} $(x_n)$ is convergent in $X$] if it does so in $(X,\tau_M).$ In notation, $x_n\to x.$ We say a sequence $(x_n)$ in $A\subset X$ is convergent in $A,$ if the same holds for $(x_n)$ in the fuzzy metric subspace $(A,M_A,*).$

\begin{theo}
\normalfont \cite{ver1} Let $(X, M, *)$ be a fuzzy metric space. A sequence $(x_n)$ in $X$ converges to $x\in X$ if and only if $M(x_n,x,t)\to 1,~\forall~t>0.$.
\end{theo}

\begin{defi}
\normalfont \cite{ver1} A sequence $(x_n)$ in a fuzzy metric space $(X, M, *)$ is said to be Cauchy if for $\epsilon\in(0,1),t>0$ there exists $k\in\mathbb N$ such that  $M(x_m,x_n,t)>1-\epsilon,~\forall~m,n\ge k.$

A fuzzy metric space is called complete if every Cauchy sequence in it converges.
\end{defi}

\begin{defi}
\normalfont A fuzzy metric space $(X,M,*)$ is said to be 

(a) compact \cite{uni} if $(X,\tau_M)$ is compact,

(b) sequentially compact \cite{sc} if a sequence in $X$ clusters to some $x\in X$ (i.e. it has a subsequence that converges to $x$).
\end{defi}

\begin{lem}\label{lemcau}
\normalfont \cite{uni} Let $(X,M,*)$ be a fuzzy metric space. If a Cauchy sequence clusters to a point $x\in X,$ then the sequence converges to $x.$
\end{lem}

We conclude this section with the notion of precompactness which provides an important tool to characterize complete fuzzy metric spaces.

\begin{defi}
\normalfont \cite{uni} A fuzzy metric space $(X,M,*)$ is said to be precompact if for $r\in(0,1),t>0$ there is a finite subset $A$ of $X$ such that $X=\bigcup_{a\in A}B(a,r,t).$
\end{defi}

\begin{lem}\label{lempre}
\normalfont \cite{uni} A fuzzy metric space $(X,M,*)$ is precompact if and only if every sequence in $X$ has a Cauchy subsequence.
\end{lem}

\begin{theo}\label{theopre}
\normalfont \cite{uni} A fuzzy metric space $(X, M, *)$ is compact if and only if it is precompact and complete.
\end{theo}

\section{\textbf{Lebesgue property for fuzzy metric spaces}}

We begin by recalling a few known facts on Lebesgue property for fuzzy metric spaces.

\begin{defi}
\normalfont \cite{uc} A fuzzy metric space $(X,M,*)$ is said to have the Lebesgue property if given an open cover $\mathcal G$ of $X$ there exist $r\in(0,1),~t>0$ such that $\{B(x,r,t):x\in X\}$ refines $\mathcal G.$ We call such fuzzy metric spaces Lebesgue.
\end{defi}

In \cite{uc}, several characterizations of Lebesgue property has been established. In what follows, we give five new characterizations of Lebesgue property for fuzzy metric spaces. In particular we show that it is sufficient to examine the uniform continuity of bounded real-valued continuous functions to characterize Lebesgue property. Among many other equivalent conditions, the following is proved in \cite{uc}: A fuzzy metric space $(X,M,*)$ if Lebesgue if and only if every real-valued continuous function on $(X, \tau_M)$ is $\mathbb R$-uniformly continuous as a mapping from $(X, M, *)$ to $\mathbb R.$

To establish the main result we require two lemmas, which are important on their own merits.

\begin{lem}\label{lemxy}
\normalfont Let $(X,M,*)$ be a fuzzy metric space. A function $f:(X,M,*)\to\mathbb R$ is $\mathbb R$-uniformly continuous if and only if given two sequences $(x_n),(y_n)$ in $X$ with $M(x_n,y_n,s)\to1,~\forall~s>0,$ we have $|f(x_n)-f(y_n)|\to0.$
\end{lem}

\noindent\textit{Proof.} \textbf{Necessity:} Let $f$ be $\mathbb R$-uniformly continuous and $(x_n),(y_n)$ be two sequences in $X$ such that $M(x_n,y_n,s)\to1,~\forall~s>0.$

Choose $\epsilon>0.$ Then there exists $r\in(0,1),s>0$ such that $$\forall~x,y\in X,~|f(x)-f(y)|<\epsilon\text{ whenever }M(x,y,s)>1-r.$$

Since $M(x_n,y_n,s)\to1,$ there exists $k\in\mathbb N$ such that $1-r<M(x_n,y_n,s)<1+r,~\forall~n\ge k.$

Thus $|f(x_n)-f(y_n)|<\epsilon,~\forall~n\ge k.$ Hence $|f(x_n)-f(y_n)|\to0.$

\textbf{Sufficiency:} Let the condition holds.

Suppose $f$ be not $\mathbb R$-uniformly continuous. Then there exists $\epsilon>0$ such that for each $n\in\mathbb N$ we can find $x_n,y_n\in X$ so that $M(x_n,y_n,\frac{1}{n})>1-\frac{1}{n}$ and $|f(x_n)-f(y_n)|\ge\epsilon.$

Choose $s>0.$ Then for any $k\in\mathbb N$ with $\frac{1}{k}<s$ we have $$1+\frac{1}{n}>M(x_n,y_n,s)\ge M(x_n,y_n,\frac{1}{n})>1-\frac{1}{n},~\forall~n\ge k.$$

Consequently $M(x_n,y_n,s)\to1,~\forall~s>0,$ a contradiction since $\lim\limits_{n\to\infty}|f(x_n)-f(y_n)|\ne0$. Hence $f$ is $\mathbb R$-uniformly continuous$.\qed$

\begin{lem}\label{lemscom}
\normalfont A fuzzy metric space $(X,M,*)$ is sequentially compact if and only if it is compact.
\end{lem}

\noindent\textit{Proof.} 
\textbf{Necessity:} Let $(X,M,*)$ be sequentially compact. Then, due to Lemma \ref{lemcau}, every Cauchy sequence in $X$ converges. Thus $X$ is complete.

In view of Theorem \ref{theopre}, it remains to show that $X$ is precompact.

If not, then there exists $r\in(0,1),t>0$ such that for no finite subset $A$ of $X$ we can have $X=\bigcup\limits_{a\in A}B(a,r,t).$ 

Choose $x_1\in X.$ Then there exists $x_2\in X\backslash~B(x_1,r,t).$ Similarly there exists $x_3\in X\backslash\bigcup\limits_{i=1}^2 B(x_i,r,t).$ Proceeding in this way we obtain a sequence $(x_n)$ of distinct elements in $X$ such that $x_n\in X\backslash\bigcup\limits_{i=1}^{n-1} B(x_i,r,t),~\forall~n\ge1.$

Clearly $M(x_m,x_n,t)\le 1-r~\forall~m\ne n$ which contradicts to the fact that $(x_n)$ clusters in $X.$ 

Hence $X$ is compact.

\textbf{Sufficiency:} Immediate from Lemma \ref{lempre} and Theorem \ref{theopre}.$\qed$

\begin{nota}
\normalfont For a fuzzy metric space $(X,M,*)$ and $A\subset X,$ we will let $B(A,r,t)$ denote $\bigcup\limits_{a\in A}B(a,r,t)$ and $M(c,A,t)$ denote $\sup\limits_{a\in A} M(c,a,t)$ where $c\in X,$ $r\in(0,1)$ and $t>0.$
\end{nota}

\begin{theo}\label{theomain}
\normalfont For a fuzzy metric space $(X,M,*),$ the following conditions are equivalent:

a) $(X,M,*)$ is Lebesgue,

b) Every real-valued continuous function on $(X,\tau_M)$ is $\mathbb R$-uniformly continuous as a mapping from $(X,M,*)$ to $\mathbb R$,

c) Every bounded real-valued continuous function on $(X,\tau_M)$ is $\mathbb R$-uniformly continuous as a mapping from $(X,M,*)$ to $\mathbb R$,

d)  Every closed and discrete subset of $(X,\tau_M)$ is fuzzy uniformly discrete in $(X, M,*),$

e) The set $X'$ of all accumulation points of $(X,\tau_M)$ is compact and for chosen $r\in(0,1),s>0,$ the set $X\backslash B(X',r,s)$ is fuzzy uniformly discrete in $(X, M,*),$

f) The set $X'$ of all accumulation points of $(X,\tau_M)$ is compact and for chosen $\delta_1\in(0,1),t_1>0$ there exist $\delta_2\in(0,1),t_2>0$ such that for each $x\in X$ with $M(x,X',t_1)\le1-\delta_1$ we have $\sup\limits_{y\ne x}M(x,y,t_2)\le1-\delta_2.$
\end{theo}

\noindent\textit{Proof.}  a$\implies$b: Follows from \cite{uc}.\\

b$\implies$c: Immediate.\\

c$\implies$d: Let $T$ be a closed and discrete subset of $(X,\tau_M).$ If possible let $T$ be not fuzzy uniformly discrete in $(X,M,*).$

We first show that $T$ is infinite. If possible let it be finite. 

Choose $s>0$ and set $p=\max\{M(x,y,s):x,y\in T,~x\ne y\}.$ Then for chosen $q\in(p,1),~M(x,y,s)\le p<q~\forall~x,y~(x\ne y)\in T.$

Setting $t=1-q\in(0,1),$ we obtain $M(x,y,s)<1-t~\forall~x,y~(x\ne y)\in T.$ But it contradicts to our assumption that $T$ is not fuzzy uniformly discrete. Hence $T$ is infinite.

Since $T$ is not fuzzy uniformly discrete, for $n>1~\exists~x_n,y_n~(x_n\ne y_n)\in T$ such that $M\left(x_n,y_n,\frac{1}{n}\right)\ge1-\frac{1}{n}.$

We now show that, for each $n>1$ we may assume $x_n,y_n\notin\{x_1,\cdots,x_{n-1},y_1,\cdots,y_{n-1}\}.$ 

If not, then there exists $v\in\mathbb N$ such that for each $m\ge v$ and $x,y~(x\ne y)\in T$ with $M\left(x,y,\frac{1}{m}\right)\ge1-\frac{1}{m},$ we have either $x\in\{x_1,\cdots,x_{v-1},y_1,\cdots,y_{v-1}\}$ or $y\in\{x_1,\cdots,x_{v-1},y_1,\cdots,y_{v-1}\}$ (since $M\left(x,y,\frac{1}{m}\right)>1-\frac{1}{m}\implies M\left(x,y,\frac{1}{v}\right)\ge M\left(x,y,\frac{1}{m}\right)>1-\frac{1}{m}\ge 1-\frac{1}{v}$).

Set $Q=\{x_1,\cdots,x_{v-1},y_1,\cdots,y_{v-1}\}.$ Then $\exists~z\in Q$ and a strictly increasing sequence $(n_k)$ in $\mathbb N$ such that $\forall~k$ there is $z_k~(\ne z)\in T$ satisfying $M(z,z_k,\frac{1}{n_k})\ge1-\frac{1}{n_k}.$ 

Choose $r\in(0,1),s>0.$ Set $n_k\in\mathbb N$ such that $\frac{1}{n_k}<\min\{r,s\}.$

Then for each $y\in B(z,\frac{1}{n_k},\frac{1}{n_k})$ we have $M(y,z,s)\ge M(y,z,\frac{1}{n_k})>1-\frac{1}{n_k}>1-r$ whence $y\in B(z,r,s).$ Consequently $B(z,\frac{1}{n_k},\frac{1}{n_k})\subset B(z,r,s).$

Since $z_k~(\ne z)\in B(z,r,s),$ $z$ is an accumulation point of $T.$

Since $T$ is closed, $z\in T.$ But it contradicts to the fact that a discrete set cannot contain its limit point. 

Thus we may assume that $x_n,y_n\notin\{x_1,\cdots,x_{n-1},y_1,\cdots,y_{n-1}\}~\forall~n\ge2.$

We now show that for $s>0,~M(x_n,y_n,s)\to1.$ 

Choose $\epsilon>0.$ Find $k\in\mathbb N$ such that $\frac{1}{k}<\min\{s,\epsilon\}.$ Then $\forall~n\ge k,$ $M(x_n,y_n,s)\ge M(x_n,y_n,\frac{1}{n})>1-\frac{1}{n}\ge 1-\frac{1}{k}>1-\epsilon.$

Thus $1-\epsilon<M(x_n,y_n,s)<1+\epsilon~\forall~n\ge k$ whence for $s>0,~M(x_n,y_n,s)\to1.$

Set $A=\{x_n:n\in\mathbb N\},~B=\{y_n:n\in\mathbb N\}.$ Clearly $A\cap B=\emptyset.$

Since $T$ is a closed, discrete subset of $(X,\tau_M)$ and $A,B\subset T$ it follows that $A,B$ are closed in $(X,\tau_M).$ 

Also $(X,\tau_M)$ is normal, being metrizable. So there exists a continuous map $f:X\to[0,1]$ such that $f(A)=\{0\},f(B)=\{1\}.$

Clearly by hypothesis $f$ is $\mathbb R$-uniformly continuous. However for all $s>0,~M(x_n,y_n,s)\to1$ even though $|f(x_n)-f(y_n)|\to1\ne 0,$ a contradiction to Lemma \ref{lemxy}.

Hence $T$ is fuzzy uniformly discrete in $X.$\\

d$\implies$e: If possible, let $X'$ be not compact. Then by Lemma \ref{lemscom}, there is a sequence $(x_n)$ in $X'$ which has no convergent subsequence.

Since $x_n\in X',~\forall~n\ge1$ there exists $y_n~(\ne x_n)\in X$ such that $y_n\in B(x_n,\frac{1}{n},\frac{1}{n}),~\forall~n\ge1.$ That is, $M(x_n,y_n,\frac{1}{n})>1-\frac{1}{n},~\forall~n\ge1.$

We first show that, $(y_n)$ has no cluster point in $X.$

Suppose not. Then $(y_n)$ has a convergent subsequence $(y_{r_n})$ such that $y_{r_n}\to l.$ 

Choose $\epsilon\in(0,1),t>0.$ Since $*$ is continuous and $1*1=1$, there exists $\delta\in(0,1)$ such that $(1-\delta)*(1-\delta)>1-\epsilon.$

Since $y_{r_n}\to l,$ there exists $k\in\mathbb N$ such that $M(y_{r_n},l,\frac{t}{2})>1-\delta,~\forall~n\ge k.$

Set $m\in\mathbb N$ such that $\frac{1}{m}<\min\{\delta,\frac{t}{2}\}.$ Then $M(x_n,y_n,\frac{t}{2})\ge M(x_n,y_n,\frac{1}{n})>1-\delta,~\forall~ n\ge m.$

Consequently $M(x_{r_n},y_{r_n},\frac{t}{2})*M(y_{r_n},l,\frac{t}{2})\ge(1-\delta)*(1-\delta),~\forall~n\ge m$ that is, $ M(x_{r_n},l,t)>1-\epsilon,~\forall~n\ge m.$ But this contradicts to the fact that $(x_n)$ has no cluster point. Thus $(y_n)$ has no cluster point. 

Therefore $B=\{x_n:n\in\mathbb N\}\cup\{y_n:n\in\mathbb N\}$ has no accumulation point in $(X,\tau_M)$ and hence is closed and discrete in $(X,\tau_M).$ So by hypothesis, $B$ is fuzzy uniformly discrete in $(X,M,*)$, a contradiction to $M(x_n,y_n,\frac{1}{n})>1-\frac{1}{n},~\forall~n\ge1.$

Hence $X'$ is compact.

Also for chosen $r\in(0,1),~s>0$ the set $X\backslash B(X',r,s)$ is closed and discrete in $(X,\tau_M)$, and hence is fuzzy uniformly discrete in $(X,M,*)$.\\

e$\implies$f: If possible, let there exist $\delta_1\in(0,1),t_1>0$ for which no such $\delta_2\in(0,1),t_2>0$ can be obtained. Then for each $n\in\mathbb N$ there exist $x_n,y_n~(x_n\ne y_n)\in X$ such that $M(x_n,X',t_1)\le1-\delta_1$ but $M(x_n,y_n,\frac{1}{n})>1-\frac{1}{n},~\forall~n\ge1.$

Since $*$ is continuous and $1*1=1,$ there exists $\delta_3\in(0,1)$ such that $x,y\in(1-\delta_3,1]\implies x*y\in(1-\delta_1,1]\cdots(1)$

Choose $k\in\mathbb N$ such that $\frac{\delta_1}{k}<\delta_3\cdots(2)$

We first show that $\exists~p\in\mathbb N$ such that $M(y_p,X',\frac{t_1}{2k})\le1-\frac{\delta_1}{2k}.$

If not then, $M(y_n,X',\frac{t_1}{2k})>1-\frac{\delta_1}{2k},~\forall~n\ge 1.$ So for $n\ge1,$ there exists $c_n\in X'$ such that $M(y_n,c_n,\frac{t_1}{2k})>1-\frac{\delta_1}{2k}~\forall~n\ge 1.$

Again $M(x_n,y_n,\frac{1}{n})>1-\frac{1}{n}~\forall~n\ge1.$

Choose $q\in\mathbb N$ such that $\frac{1}{q}<\min\{\frac{\delta_1}{2k},\frac{t_1}{2k}\}.$ 

Then $\forall~n\ge q,$ $M(x_n,y_n,\frac{t_1}{2k})\ge M(x_n,y_n,\frac{1}{q})\ge M(x_n,y_n,\frac{1}{n})>1-\frac{1}{n}\ge1-\frac{1}{q}>1-\frac{\delta_1}{2k}.$

Thus in view of $(1)$ and $(2),$ $M(y_n,c_n,\frac{t_1}{2k})*M(x_n,y_n,\frac{t_1}{2k})\ge(1-\frac{\delta_1}{2k})*(1-\frac{\delta_1}{2k})>1-\delta_1,~\forall~n\ge q$ whence $M(x_n,c_n,t_1)>1-\delta_1,~\forall~n\ge q.$

Consequently $M(x_n,X',t_1)>1-\delta_1,~\forall~n\ge q,$ a contradiction.

Hence for some $p\in\mathbb N,$ we have $M(y_p,X',\frac{t_1}{2k})\le1-\frac{\delta_1}{2k}.$

Similarly, there exists $p'~(>p)\in\mathbb N$ such that $M(y_{p'},X',\frac{t_1}{2k})\le1-\frac{\delta_1}{2k}.$

Continuing in this way we may pass $(y_n)$ to a subsequence such that $M(y_n,X',\frac{t_1}{2k})\le1-\frac{\delta_1}{2k}~\forall~n\ge1.$

Since $\forall~n\ge1\text{ and }c\in X',$ $M(y_n,c,\frac{t_1}{2k})\le1-\frac{\delta_1}{2k}$ so we have, $y_n\notin B(X',\frac{\delta_1}{2k},\frac{t_1}{2k}),~\forall~n\ge1.$

Again $\forall~n\ge1\text{ and }c\in X',$ $M(x_n,c,t_1)\le1-\delta_1\implies x_n\notin B(X',\frac{\delta_1}{2k},\frac{t_1}{2k}),~\forall~n\ge1.$

Thus $x_n,y_n\in X\backslash B(X',\frac{\delta_1}{2k},\frac{t_1}{2k})~\forall~n\ge1.$ 

Since by hypothesis, $X\backslash B(X',\frac{\delta_1}{2k},\frac{t_1}{2k})$ is uniformly discrete, there exist $r\in(0,1),s>0$ such that $M(x_n,y_n,s)<1-r~\forall~n\ge1.$

Choose $v\in\mathbb N$ such that $\frac{1}{v}<\min\{r,s\}.$ Then $M(x_v,y_v,\frac{1}{v})\le M(x_v,y_v,s)<1-r<1-\frac{1}{v},$ a contradiction to the fact that $M(x_n,y_n,\frac{1}{n})<1-\frac{1}{n}~\forall~n\ge1.$

Hence the result follows.\\

f$\implies$a: Let $\{O_\lambda:\lambda:\Lambda\}$ be an open cover of $X.$ We first show that there is $\delta\in(0,1),t>0$ such that for $x\in X'~\exists~\lambda_x\in\Lambda$ such that $B(x,\delta,t)\subset O_{\lambda_x}.$

If not then, for each $\delta=t=\frac{1}{n}~(n\in\mathbb N\backslash\{1\})$ there is $x_n\in X'$ such that $B(x_n,\frac{1}{n},\frac{1}{n})\not\subset O_{\lambda}~\forall~\lambda\in\Lambda.$

Since $X'$ is compact, by Lemma \ref{lemscom}, $(x_n)$ has a subsequence $(x_{r_n})$ converging to some point $w\in X'.$ Choose $\lambda_1\in\Lambda$ such that $w\in O_{\lambda_1}.$

Since $O_{\lambda_1}$ is open, there exist $\delta'\in(0,1),t>0$ such that $B(x,\delta',t')\subset O_{\lambda_1}.$

Choose $q>1$ such that $\frac{1}{r_q}<\min\{\delta',t'\}.$ Then $B(w,\frac{1}{r_q},\frac{1}{r_q})\subset O_{\lambda_1}$ since $y\in B(w,\frac{1}{r_q},\frac{1}{r_q})\implies M(y,w,t')\ge M(y,w,\frac{1}{r_q})>1-\frac{1}{r_q}>1-\delta'\implies y\in B(w,\delta',t')\subset O_{\lambda_1}.$

Since $*$ is continuous and $1*1=1,~\exists~\delta_1\in(0,1)$ such that $x',y'\in(1-\delta_1,1]\implies x'*y'\in(1-\frac{1}{2r_q},1].$

Choose $p~(>q)\in\mathbb N$ such that $1-\frac{1}{2r_p}>1-\delta_1\cdots(1).$ 

Since $x_{2r_n}\to w,$ there is $s~(\ge p)\in\mathbb N$ such that $x_{2r_n}\in B(w,\frac{1}{2r_p},\frac{1}{2r_p}),~\forall~n\ge s.$ That is, $M(x_{2r_n},w,\frac{1}{2r_p})>1-\frac{1}{2r_p},~\forall~n\ge s\cdots(2)$

Choose $v\ge s.$ 

Let $y\in B(x_{2r_v},\frac{1}{{2r_v}},\frac{1}{{2r_v}}).$ Then $M(x_{2r_v},y,\frac{1}{2r_p})\ge M(x_{2r_v},y,\frac{1}{2r_v})>1-\frac{1}{2r_v}\ge1-\frac{1}{2r_p}\implies M(x_{2r_v},y,\frac{1}{2r_p})>1-\frac{1}{2r_p}\cdots(3)$

Again using $(2)$ we obtain, $M(x_{2r_v},w,\frac{1}{2r_p})>1-\frac{1}{2r_p}\cdots(4)$

Thus from $(3)$ and $(4)$ we have,

$M(x_{2r_v},y,\frac{1}{2r_p})*M(x_{2r_v},w,\frac{1}{2r_p})\ge(1-\frac{1}{2r_p})*(1-\frac{1}{2r_p})\implies M(w,y,\frac{1}{r_p})>1-\frac{1}{2r_q}>1-\frac{1}{r_q}$ [using $(1)$]. So, $y\in B(w,\frac{1}{r_q},\frac{1}{r_q}).$ 

Thus $B(x_{2r_v},\frac{1}{2r_v},\frac{1}{2r_v})\subset B(w,\frac{1}{r_q},\frac{1}{r_q})\subset O_{\lambda_1},$ a contradiction.

Hence there exist $\delta\in(0,1),t>0$ such that given $x\in X',$ we have $B(x,\delta,t)\subset O_{\lambda_x}$ for some $\lambda_x\in\Lambda.$ 

Choose $\rho\in(0,1)$ such that $x',y'\in(1-\rho,1]\implies x'*y'\in(1-\delta,1].$

Let $x\in B(X',\frac{\rho}{2},\frac{t}{2}).$ Then $x\in B(x',\frac{\rho}{2},\frac{t}{2})$ for some $x'\in X'.$

Choose $y\in B(x,\frac{\rho}{2},\frac{t}{2}).$ Then $M(x,y,\frac{t}{2})>1-\frac{\rho}{2}.$ Also $M(x,x',\frac{t}{2})>1-\frac{\rho}{2}.$ Thus $M(y,x',t)\ge M(x,y,\frac{t}{2})*M(x,x',\frac{t}{2})\ge(1-\frac{\rho}{2})*(1-\frac{\rho}{2})>1-\delta\implies y\in B(x',\delta,t).$ So $B(x,\frac{\rho}{2},\frac{t}{2})\subset B(x',\delta,t)\subset O_{\lambda_{x'}}.$

Note that $x\in X\backslash B(X',\frac{\rho}{2},\frac{t}{2})\implies M(x,x',\frac{t}{2})\le 1-\frac{\rho}{2},~\forall~x'\in X'\implies M(x,X',\frac{t}{2})\le1-\frac{\rho}{2}.$ So by hypothesis $\exists~\alpha\in(0,1),\beta>0$ such that $\sup\limits_{y\ne x} M(x,y,\beta)\le1-\alpha,~\forall~x\in X\backslash B(X',\frac{\rho}{2},\frac{t}{2}).$

Clearly $B(x,\alpha,\beta)=\{x\},~\forall~x\in X\backslash B(X',\frac{\rho}{2},\frac{t}{2}).$ Thus for $x\in X\backslash B(X',\frac{\rho}{2},\frac{t}{2}),$ $B(x,\alpha,\beta)\subset O_\lambda$ for some $\lambda\in\Lambda.$

Set $\gamma=\min\{\frac{\rho}{2},\alpha\},\eta=\min\{\frac{t}{2},\beta\}.$ Then it is easy to see that, for $x\in X,$ $B(x,\gamma,\eta)\subset O_\lambda$ for some $\lambda\in\Lambda.$

Hence the result follows.\\

The following example illustrates that a closed and discrete subset of a metric space may not be fuzzy uniformly discrete in the induced fuzzy metric space.

\begin{exm}\label{exm1}
\normalfont Suppose $Y=(0,1]$ be endowed with the usual metric $\rho$ on it. Set $A=\left\{\frac{1}{n}:n\in\mathbb N\right\}.$ Then $A$ is a closed and discrete subset of $(Y,\rho).$

However $A$ is not fuzzy uniformly discrete in $(Y,M_\rho,*).$ For otherwise, there exists $r\in(0,1),~t>0$ such that $$\frac{t}{t+|x-y|}<1-r$$ for all $x,y~(x\ne y)\in A.$ This is a contradiction since $|x-y|$ can be made arbitrarily small for $x,y\in A$. Thus $A$ is not fuzzy uniformly discrete in $(Y,M_\rho,*).$

Hence, in view of Theorem \ref{theomain}, $(Y,M_\rho,*)$ is not Lebesgue.
\end{exm}


\begin{thebibliography}{99}

\bibitem{1} Deng Zi-ke, Fuzzy pseudo-metric spaces, J. Math. Anal. Appl. 86 (1982) 74-95.

\bibitem{2} M.A. Erceg, Metric spaces in fuzzy set theory, J. Math. Anal. Appl. 69 (1979) 205-230.

\bibitem{ver1} A. George, P. V. Veeramani, On some results in fuzzy metric spaces, Fuzzy Sets and Systems 64 (1994) 395-399.

\bibitem{ver2}  A. George, P.V. Veeramani, On some results of analysis for fuzzy metric spaces, Fuzzy Sets and Systems 90 (1997) 365-368.

\bibitem{cont} A. George, P.V. Veeramani, Some theorems in fuzzy metric spaces, J. Fuzzy Math. 3 (1995) 933-940.

\bibitem{3} O. Kaleva, S. Seikkala, On fuzzy metric spaces, Fuzzy Sets and Systems 12 (1984) 215-229.

\bibitem{km} O. Kramosil, J. Michalek, Fuzzy metric and statistical metric spaces, Kybernetica 11 (1975) 326-334.

\bibitem{sub} L. D.R. Ko\v{c}inac, Selection properties in fuzzy metric spaces, Filomat 26 (2012), 305-312.

\bibitem{lim} V. Gregori, A. L\'{o}pez-Crevill\'{e}n, S. Morillas, A. Sapena, On convergence in fuzzy metric spaces, Top. App. 156 (2009) 3002-3006.

\bibitem{uni} V. Gregori, S. Romaguera, Some properties of fuzzy metric spaces, Fuzzy Sets and Systems 115 (2000) 485-489.

\bibitem{uc} V. Gregori, S. Romaguera, A. Sapena, Uniform continuity in fuzzy metric spaces, Rend. Ist. Mat. Univ. Trieste 32 Suppl. 2 (2001) 81-88.

\bibitem{tn} B. Schweizer, A. Sklar, Statistical metric spaces, Pacific J. Math. 10 (1960) 314-334.

\bibitem{dis}  P. Veeramani, Best approximation in fuzzy metric spaces, J. Fuzzy Math. 9 (2001) 75-80.

\bibitem{fd} C. Li, Y. Zhang, On precompactness of the Hausdorff fuzzy metric on closed sets, Journal of Comp. Anal. and App. 24 (2018) 343-353.

\bibitem{sc} K. P. R. Rao, K. R. K. Rao, T. Ranga Rao, Common fixed point theorems in sequentially compact fuzzy metric spaces, Int. Math. Forum 2 (2007) 2543 - 2549.

\bibitem{will} S. Willard, General Topology. Reading, Mass.: Addison-Wesley Publishing, 1970.

\end{thebibliography}
\end{document}